\documentclass[11pt]{article}
\usepackage{times}
\usepackage{amssymb,amsmath}
\usepackage{epsf}
\usepackage{times,bm,mathrsfs}  
\usepackage{amsmath}
\usepackage{amssymb}
\usepackage{amsfonts}
\usepackage{mathrsfs}
\usepackage{bbm}
\usepackage{color}
\usepackage{graphicx}
\usepackage{amscd}
\usepackage{epsfig}

\definecolor{Red}{rgb}{1,0,0}
\definecolor{Blue}{rgb}{0,0,1}
\definecolor{Olive}{rgb}{0.41,0.55,0.13}
\definecolor{Green}{rgb}{0,1,0}
\definecolor{MGreen}{rgb}{0,0.8,0}
\definecolor{DGreen}{rgb}{0,0.55,0}
\definecolor{Yellow}{rgb}{1,1,0}
\definecolor{Cyan}{rgb}{0,1,1}
\definecolor{Magenta}{rgb}{1,0,1}
\definecolor{Orange}{rgb}{1,.5,0}
\definecolor{Violet}{rgb}{.5,0,.5}
\definecolor{Purple}{rgb}{.75,0,.25}
\definecolor{Brown}{rgb}{.75,.5,.25}
\definecolor{Grey}{rgb}{.5,.5,.5}
\definecolor{Black}{rgb}{0,0,0}

\def\path{{\tt path}}

\newcommand{\eps}{\varepsilon}

\newcommand{\bdm}{\begin{displaymath}}
\newcommand{\edm}{\end{displaymath}}
\newcommand{\bea}{\begin{eqnarray*}}
\newcommand{\eea}{\end{eqnarray*}}
\newcommand{\bean}{\begin{eqnarray}}
\newcommand{\eean}{\end{eqnarray}}

\newcommand{\polylog}{\mathrm{polylog}}

\newtheorem{theorem}{Theorem}[section]
\newtheorem{lemma}[theorem]{Lemma}
\newtheorem{proposition}[theorem]{Proposition}
\newtheorem{corollary}[theorem]{Corollary}

\newtheorem{definition}[theorem]{Definition}
\newtheorem{remark}[theorem]{Remark}
\newtheorem{conjecture}[theorem]{Conjecture}
\newtheorem{algorithm}[theorem]{Algorithm}

\newenvironment{proof}{\noindent{\textbf{Proof:}}}{$\blacksquare$\vskip\belowdisplayskip}

\author{
DRAFT - DO NOT CIRCULATE
}
\title{Rapid Mixing of Gibbs Sampling on Graphs that are Sparse on Average
}

\author{Elchanan Mossel\thanks{Email: mossel@stat.berkeley.edu. Dept. of Statistics, U.C. Berkeley. Supported by an Alfred Sloan fellowship
  in Mathematics and by NSF grants DMS-0528488,
  DMS-0548249 (CAREER) by  DOD ONR grant N0014-07-1-05-06.} \and
  Allan Sly \thanks{{Email: sly@stat.berkeley.edu Dept. of Statistics, U.C. Berkeley.}}}
\begin{document}

\maketitle

\thispagestyle{empty}

\begin{abstract}

Gibbs sampling also known as Glauber dynamics
is a popular technique for sampling high
dimensional distributions defined on graphs. Of special interest is the
behavior of Gibbs sampling on the Erd\H{o}s-R\'enyi random graph $G(n,d/n)$, where
each edge is chosen independently with probability $d/n$ and $d$ is fixed.
While the average degree in $G(n,d/n)$ is $d(1-o(1))$,
it contains many nodes of degree of order $\log n / \log \log n$.

The existence of nodes of almost logarithmic degrees implies that for many
natural distributions defined on $G(n,p)$ such as uniform coloring
(with a constant number of colors) or the Ising model at any fixed
inverse temperature $\beta$, the mixing time of Gibbs sampling is at
least $n^{1 + \Omega(1 / \log \log n)}$.
Recall that the Ising model with inverse temperature
$\beta$ defined on a
graph $G=(V,E)$ is the distribution over $\{ \pm \}^{V}$ given by
$P(\sigma) = \frac{1}{Z} \exp(\beta \sum_{(v,u) \in E} \sigma(v) \sigma(u))$.
High degree nodes
pose a technical challenge in proving polynomial time mixing of the dynamics
for many models including the Ising model and coloring. Almost all known
sufficient conditions in terms of $\beta$ or number of colors needed
for rapid mixing of Gibbs samplers are stated in terms
of the maximum degree of the underlying graph.

In this work we show that for every $d < \infty$
and the Ising model defined on
$G(n,d/n)$, there exists a $\beta_d > 0$, such that for all
$\beta < \beta_d$ with probability going to $1$ as $n \to \infty$, the mixing
time of the dynamics on $G(n,d/n)$ is polynomial in $n$.
Our results are the first polynomial time mixing results
proven for a natural model on $G(n,d/n)$ for $d > 1$
where the parameters of the model do not depend on $n$. They also provide a rare example where one can prove a polynomial time mixing of Gibbs sampler
in a situation where the actual mixing time is slower than
$n \polylog(n)$.
Our proof exploits in novel ways the local treelike structure of Erd\H{o}s-R\'enyi
random graphs, comparison and block dynamics arguments and a recent
result of Weitz. 

Our results extend to much more general families of graphs which are
sparse in some average sense and to much more general interactions.
In particular, they apply to any graph for which every vertex $v$ of the graph
has a neighborhood $N(v)$ of radius $O(\log n)$
in which the induced sub-graph is a tree
union at most $O(\log n)$ edges and where for each simple path in $N(v)$ the
sum of the vertex degrees along the path is $O(\log n)$. Moreover, our result apply also in the case of arbitrary external fields and provide the first FPRAS
for sampling the Ising distribution in this case. We finally present a non 
Markov Chain algorithm for sampling the distribution which is effective for a 
wider range of parameters. In particular, for $G(n,d/n)$ it applies 
for all external fields and $\beta < \beta_d$, 
where $d \tanh(\beta_d) = 1$ is the critical 
point for decay of correlation for the Ising model on $G(n,d/n)$.

\end{abstract}

\bigskip


\clearpage

\section{Introduction}
Efficient approximate sampling from Gibbs distributions is a central
challenge of randomized algorithms. Examples include sampling from
the uniform distribution over independent sets of a
graph~\cite{Weitz:06,Vigoda:01,DyFrJe:02,DyerGreenhill:97}, sampling
from the uniform distribution of perfect matchings in a
graph~\cite{JeSiVi:04}, or sampling from the uniform distribution of
colorings~\cite{GoMaPa:04,DFHV:04,DFFV:06} of a graph. A natural
family of approximate sampling techniques is given by Gibbs
samplers, also known as Glauber dynamics. These are reversible
Markov chains that have the desired distribution as their stationary
distribution and where at each step the status of one vertex is
updated. It is typically easy to establish that the chains
will eventually converge to the desired distribution.\\

Studying the convergence rate of the dynamics is interesting from
both the theoretical computer science the statistical physics
perspectives. Approximate convergence in time polynomial in the size of the system, sometimes
called {\em rapid mixing}, is essential in computer science
applications. The convergence rate is also of natural interest in physics where the dynamical properties of such distributions are
extensively studied, see e.g.~\cite{Martinelli:99}. Much recent work
has been devoted to determining sufficient and necessary conditions
for rapid convergence of Gibbs samplers. A common feature to most of
this
work~\cite{Weitz:06,Vigoda:01,DyFrJe:02,DyerGreenhill:97,GoMaPa:04,DFHV:04,KeMoPe:01,MaSiWe:04}
is that the conditions for convergence are stated in terms of the
maximal degree of the underlying graph. In particular, these results
do not allow for the analysis of the mixing rate of Gibbs samplers
on the Erd\H{o}s-R\'enyi random graph, which is sparse on average,
but has rare denser sub-graphs.
 Recent work has been directed at showing
how to relax statements so that they do not involve maximal
degrees~\cite{DFFV:06,Hayes:06}, but the results are not strong
enough to imply rapid mixing of Gibbs sampling for the Ising model
on $G(n,d/n)$ for $d
> 1$ and any $\beta > 0$ or for sampling uniform colorings from
$G(n,d/n)$ for $d > 1$ and $1000d$ colors. The second challenge
is presented as the major open problem of~\cite{DFFV:06}.\\

In this paper we give the first rapid convergence result of Gibbs samplers
for the Ising model on Erd\H{o}s-R\'enyi random graphs in terms of the average degree and $\beta$ only.
Our results hold for the Ising model allowing different interactions and
arbitrary external fields.
We note that there is an FPRAS that samples from the Ising model on any graph~\cite{JerrumSinclair:93} as long as all the interactions are positive and the external field is the same for all vertices. However, these results do not provide
a FPRAS in the case where different nodes have different external fields as we do here.

Our results are further extended to much more general families of
graphs that are ``tree-like'' and ``sparse on average''. These are
graph where every vertex has a radius $O(\log n)$ neighborhood which
is a tree with at most $O(\log n)$ edges added and where for each
simple path in the neighborhood, the sum of degrees along the path
is $O(\log n)$. An important open problem~\cite{DFFV:06} is to
establish similar conditions for other models defined on graphs,
such as the uniform distribution over colorings.

Below we define the Ising model and Gibbs samplers and state our main result.
Some related work and a sketch of the proof are also given as the introduction.
Section~\ref{sec:proofs} gives a more detailed proof though we have not tried
to optimize any of the parameters in proofs below.

\subsection{The Ising Model}

The Ising model is perhaps the simplest model defined on graphs.
This model defines a distribution on labelings of the vertices of
the graph by $+$ and $-$. The Ising model has various natural
generalizations including the uniform distribution over colorings.
The Ising model with varying parameters is of use in a variety of
areas of machine learning, most notably in vision, see
e.g.~\cite{GemanGeman:84}.

\begin{definition}
The (homogeneous)  Ising model on a (weighted) graph $G$ with inverse temperature
$\beta$ is a distribution on configurations $\{\pm\}^V$ such that
\begin{equation} \label{eq:def_ising}
P(\sigma)=\frac1{Z(\beta)}\exp(\beta\sum_{\{v,u\}\in
E}\sigma(v)\sigma(u))
\end{equation}
where $Z(\beta)$ is a normalizing constant.

More generally, we will be interested in (inhomogeneous) Ising models defined by:
\begin{equation} \label{eq:def_ising_general}
P(\sigma)=\frac1{Z(\beta)}\exp(\sum_{\{v,u\}\in
E} \beta_{u,v} \sigma(v)\sigma(u) + \sum_v h_v \sigma(v)),
\end{equation}
where $h_v$ are arbitrary and where $\beta_{u,v} \geq 0$ for all $u$ and $v$.
In the more general case we will write $\beta = \max_{u,v} \beta_{u,v}$.
\end{definition}

\subsection{Gibbs Sampling}

The Gibbs sampler is a Markov chain on configurations where a
configuration $\sigma$ is updated by choosing a vertex $v$ uniformly
at random and assigning it a spin according to the Gibbs
distribution conditional on the spins on $G-\{v\}$.
\begin{definition}
Given a graph $G=(V,E)$ and an inverse temperature $\beta$, the
Gibbs sampler is the discrete time Markov chain on $\{\pm\}^V$ where given
the current configuration $\sigma$ the next configuration $\sigma'$ is obtained
by choosing a vertex $v$ in $V$ uniformly at random and
\begin{itemize}
\item
Letting $\sigma'(w) = \sigma(w)$ for all $w \neq v$.
\item
$\sigma'(v)$ is
assigned the spin $+$ with probability
\[
\frac{\exp(h_v + \sum_{u:(v,u)\in
E}\beta_{u,v} \sigma(u))}{\exp(h_v + \sum_{u:(v,u)\in
E}\beta_{u,v} \sigma(u))+\exp(-h_v-\sum_{u:(v,u)\in E} \beta_{u,v} \sigma(u))}.
\]
\end{itemize}
\end{definition}

We will be interested in the time it takes the dynamics to get close
to the distributions~(\ref{eq:def_ising})
and~(\ref{eq:def_ising_general}). The {\em mixing time }
$\tau_{mix}$ of the chain is defined as the number of steps needed
in order to guarantee that the chain, starting from an arbitrary
state, is within total variation distance $1/2e$ from the
stationary distribution. We will bound the mixing time by the relaxation time defined below.\\

It is well known that Gibbs sampling is a reversible Markov chain
with stationary distribution $P$.  Let $1=\lambda_1 > \lambda_2 \geq
\ldots \geq \lambda_m \geq -1$ denote the eigenvalues of the
transition matrix of Gibbs sampling.  The {\em spectral gap} is
denoted by $\max\{1-\lambda_2,1-|\lambda_m|\}$ and the {\em
relaxation time} $\tau$ is the inverse of the spectral gap.  The
relaxation time can be given in terms of the Dirichlet form of the
Markov chain by the equation
\begin{equation}\label{eq_relax_defn}
\tau=\sup\left\{\frac{2\sum_\sigma P(\sigma)
(f(\sigma))^2}{\sum_{\sigma\neq\tau} Q(\sigma,\tau)
(f(\sigma)-f(\tau))^2} : \sum_\sigma P(\sigma)
f(\sigma) \neq 0 \right \}
\end{equation}
where $f:\{\pm\}^V\rightarrow \mathbb{R}$ is any function on configurations,
$Q(\sigma,\tau)=P(\sigma)P(\sigma \rightarrow \tau)$ and $P(\sigma
\rightarrow \tau)$ is transition probability from $\sigma$ to
$\tau$. We use the result that for reversible Markov chains the
relaxation time satisfies
\begin{equation} \label{eq:tau_and_spectral}
\tau \leq \tau_{mix}\leq  \tau\left(1+\frac12 \log(\min_\sigma
P(\sigma))^{-1}\right)
\end{equation}
where $\tau_{mix}$ is the mixing time (see e.g. \cite{AldousFill:u})
and so by bounding the relaxation time we can bound the mixing time
up to a polynomial factor.\\

While our results are given for the discrete time Gibbs Sampler described above, it will at times be convenient to consider the continuos time version of the model.  Here sites are updated at rate 1 by independent Poisson clocks.  The two chains are closely related, the relaxation time of the continuous time Markov chain is $n$ times the relaxation time of the discrete chain (see e.g. \cite{AldousFill:u}).

For our proofs it will be useful to use the notion of {\em block
dynamics}. The Gibbs sampler can be generalized to update blocks of
vertices rather than individual vertices.  For blocks
$V_1,V_2,\ldots,V_k\subset V$ with $V=\cup_i V_i$ the block dynamics
of the Gibbs sampler updates a configuration $\sigma$ by choosing a
block $V_i$ uniformly at random and assigning the spins in $V_i$
according to the Gibbs distribution conditional on the spins on
$G-\{V_i\}$.  There is also a continuous analog in which the blocks each update at rate 1.  In continuous time, the relaxation time of the Gibbs sampler can be
given in terms of the relaxation time of the block dynamics and the
relaxation times of the Gibbs sampler on the blocks.
\begin{proposition} \label{prop:martinelli}
In continuous time if $\tau_{block}$ is the relaxation time of the block dynamics and
$\tau_{i}$ is the maximum the relaxation time on $V_i$ given any
boundary condition from $G-\{V_i\}$ then by Proposition 3.4 of
\cite{Martinelli:99}
\begin{equation}\label{eq_relax_block}
\tau\leq \tau_{block}(\max_i \tau_{i})\max_{v\in V} \{\# j: v\in
V_j\}.
\end{equation}
\end{proposition}

\subsubsection{Monotone Coupling}\label{s:monotoneCoupling}
For two configurations $X,Y\in\{ - , + \}^V$ we let $X \succcurlyeq Y$ denote that $X$ is greater than or equal to $Y$ pointwise.  When all the interactions $\beta_{ij}$ are positive, it is well known that the Ising model is a monotone system under this partial ordering, that is if $X \succcurlyeq Y$ then,
\[
P\left(\sigma_v = + | \sigma_{V \setminus \{ v \}} = X_{V \setminus \{ v \} }\right) \geq P\left(\sigma_v = + | \sigma_{V \setminus \{ v \}} =Y_{V \setminus \{ v \} } \right).
\]

As it is a monotne system there exists a coupling of Markov chains $\{X^x_t\}_{x\in \{ -,+\}^V } $ such that marginally each has the law of the Gibbs Sampler with starting configurations $X_0^x = X$ and further that if $x \succcurlyeq y$ then for all $t$, $X_t^x \succcurlyeq X_t^y$.  This is referred to as the monotone coupling and can be constructed as follows: let $v_1,\ldots$ be a random sequence of vertices updated by the Gibbs Sampler and associate with them iid random variables $U_1,\ldots$ distributed as $U[0,1]$ which determine how the site is updated.  At the $i$th update the site $v_i$ is updated to $+$ if $$U_i \leq \frac{\exp(h_v + \sum_{u:(v,u)\in
E}\beta_{u,v} \sigma(u))}{\exp(h_v + \sum_{u:(v,u)\in
E}\beta_{u,v} \sigma(u))+\exp(-h_v-\sum_{u:(v,u)\in E} \beta_{u,v} \sigma(u))}
$$
and to $-$ otherwise.  It is well known that such transitions preserve the partial ordering  which guarantees that if $x \succcurlyeq y$ then $X_t^x  \succcurlyeq   X^y_t$ by the monotonicity of the system.  In particular this implies that it is enough to bounded the time taken to couple from the all $+$ and all $-$ starting configurations.

\subsection{Erd\H{o}s-R\'enyi Random Graphs and Other Models of graphs}

The Erd\H{o}s-R\'enyi random graph $G(n,p)$, is the graph with $n$ vertices
$V$ and random edges $E$ where each potential edge $(u,v) \in V \times V$
is chosen independently with probability $p$. We take
$p=d/n$ where $d \geq 1$ is fixed. In the case $d < 1$, it is well known that
with high probability all components of $G(n,p)$ are of logarithmic size which implies immediately that the dynamics mix in polynomial time for all $\beta$.

For a vertex $v$ in $G(n,d/n)$ let $V(v,l)=\{u\in G: d(u,v)\leq
l\}$, the set of vertices within distance $l$ of $v$, let
$S(v,l)=\{u\in G:d(u,v)=l\}$, let $E(v,l)=\{ (u,w)\in G: u,w \in
V(v,l) \}$ and let $B(v,l)$ be the graph ($V(v,l),E(v,l))$.

Our results only require some simple features of the neighborhoods of all vertices in the graph.
\begin{definition}
Let $G = (V,E)$ be a graph and $v$ a vertex in $G$.
Let $t(G)$ denote the {\em tree access} of $G$, i.e.,
\[
t(G) = |E| - |V| + 1.
\]
We call a path $v_1,v_2,\ldots$ {\em self avoiding} if
for all $i\neq j$ it holds that $v_i \neq v_{j}$. We let the {\em maximal
path density} $m$ be defined by
\[
m(G,v,l) =
\max_{\Gamma} \sum_{u \in \Gamma} d_u
\]
where the maximum is taken
over all self-avoiding paths $\Gamma$ starting at $v$ with length
at most $l$ and $d_u$ is the degree of node $u$. We write $t(v,l)$
for $t(B(v,l))$ and $m(v,l)$ for $m(B(v,l),v,l)$.
\end{definition}

\subsection{Our Results}

Throughout we will be using the term \emph{with high probability} to mean with probability $1-o(1)$ as $n$ goes to $\infty$.

\begin{theorem}\label{thm_mix_time_random_graph}
Let $G$ be a random graph distributed as $G(n,d/n)$.  When
\[
\tanh(\beta)< \frac{1}{e^2 d},
\]
there exists constant a $C=C(d)$
such that the mixing time of the Glauber dynamics is $O(n^{C})$ with
high probability (probability $1-o(1)$) over the graph as $n$ goes to $\infty$.
The result holds for the homogeneous model~(\ref{eq:def_ising})
and for the inhomogeneous model~(\ref{eq:def_ising_general})
provided $|h_v| \leq 100 \beta n$ for all $v$.
\end{theorem}
Note in the theorem above the $O(\cdot)$ bound depends on $\beta$. 
It may be viewed as a special case of the following more general result.

\begin{theorem}\label{thm_mix_general_graph}
Let $G=(V,E)$ be any graph on $n$ vertices satisfying the following properties.
There exist $a >0, 0 < b < \infty$ and $0 < c < \infty$ such that for all $v \in V$ it holds that
\[
t(v,a \log n)\leq b \log n, \quad
m(v,a \log n)\leq c \log n.
\]
Then if
\[
\tanh(\beta)< \frac{a}{e^{1/a} (c-a)},
\]
there exists constant a $C=C(a,b,c,\beta)$
such that the mixing time of the Glauber dynamics is $O(n^{C})$.
The result holds for the homogeneous model~(\ref{eq:def_ising})
and for the inhomogeneous model~(\ref{eq:def_ising_general})
provided $|h_v| \leq 100 \beta n$ for all $v$.
\end{theorem}

\begin{remark} \label{remark:h}
The condition that $|h_v| \leq 100 \beta n$ for all $v$ will be
needed in the proof of the result in the general
case~(\ref{eq:def_ising_general}). However, we note that given
Theorem~\ref{thm_mix_general_graph} as a black box, it is easy to
extend the result and provide an efficient sampling algorithm in the
general case without any bounds on the $h_v$. In the case where some
of the vertices $v$ satisfy $|h_v| \geq 10 \beta n$, it is easy to
see that the target distribution satisfies except with exponentially
small probability that $\sigma_v = +$ for all $v$ with $h_v > 10
\beta n$ and $\sigma_v = -$ for all $v$ with $h_v < -10 \beta n$.
Thus we may set $\sigma_v = +$ when $h_v > 10 \beta n$ and $\sigma_v
= -$ when $h_v < 10 \beta n$ and consider the dynamics where these
values are fixed. Doing so will effectively restrict the dynamics
to the graph spanned by the remaining vertices and will modify the
values of $h_v$ for the remaining vertices; however, it is easy to
see that all remaining vertices will have $|h_v| \leq 100 \beta n$.
It is also easy to verify that if the original graph satisfied the
hypothesis of Theorem~\ref{thm_mix_general_graph} then so does the
restricted one. Therefore we obtain an efficient sampling procedure
for the desired distribution.
\end{remark}

\subsection{Related Work and Open Problems}

Much work has been focused on the problem of understanding the mixing time of the Ising model in various contexts.  In a series of results \cite{Holley:85,AizHol:87, Zegarlinski:90} culminating in \cite{StrookZegarlinski:92} it was shown that the Gibbs sampler on integer lattice mixes rapidly when the model has the strong spatial mixing property.  In $\mathbb{Z}^2$ strong spatial mixing, and therefore rapid mixing, holds in the entire uniqueness regime (see e.g. \cite{MartinelliOlieri:94a}).  On the regular tree the mixing time is always polynomial but is only $O(n\log n)$ up to the threshold for extremity \cite{BeKeMoPe:05}.  For completely general graphs the best known results are given by the Dobrushin condition which establishes rapid mixing when $d\tanh(\beta)<1$ where $d$ is the maximum degree.

Most results for mixing rates of Gibbs samplers are stated in terms
of the maximal degree. For example many results have focused on  sampling uniform colorings, the
result are of the form: for every graph where all degrees are at
most $d$ if the number of colors $q$ satisfies $q \geq q(d)$ then
Gibbs sampling is rapidly
mixing~\cite{Weitz:06,Vigoda:01,DyFrJe:02,DyerGreenhill:97,GoMaPa:04,DFHV:04,KeMoPe:01,MaSiWe:04}.
For example, Jerrum \cite{Jerrum:95} showed that one can take
$q(d) = 2d$. The novelty of the result presented
here is that it allows for the study of graphs where the average degree is
small while some degrees may be large.

Previous attempts at studying this problem, with bounded average degree but some large degrees, for sampling uniform
colorings yielded weaker results. In~\cite{DFFV:06} it is shown that
Gibbs sampling rapidly mixes on $G(n,d/n)$ if $q = \Omega_d((\log
n)^{\alpha})$ where $\alpha < 1$ and that a variant of the algorithm
rapidly mixes if $q \geq \Omega_d(\log \log n / \log \log \log n)$.
Indeed the main open problem of~\cite{DFFV:06} is to determine if
one can take $q$ to be a function of $d$ only. Our results here
provide a positive answer to the analogous question for the Ising
model. We further note that other results where the conditions on
degree are
relaxed~\cite{Hayes:06} do not apply in our setting.\\

The following propositions, which are easy and well known, establish that for $d > 1$ and large $\beta$ the mixing time is exponential in $n$
and that for all $d >0$ and
$\beta > 0$ the mixing time is more than $n \polylog(n)$.

\begin{proposition}
If $d >0$ and $\beta > 0$ then with high probability
the mixing time of the dynamics on $G(n,d/n)$ is at least
$n^{1+\Omega(1/\log \log n)}$.
\end{proposition}

\begin{proof}
The proof follows from the fact that $G(n,d/n)$ contains an isolated
star with $s = \Omega(\log n / \log \log n)$ vertices with high
probability and that the mixing time of the star is $s
\exp(\Omega(s))$. Since the star is updated with frequency $s/n$, it
follows that the mixing time is at least
\[
(n/s) s \exp(\Omega(s)) = n \exp(\Omega(s)) = n^{1+\Omega(1/\log \log n)}.
\]
\end{proof}

\begin{proposition}
If $d > 1$ then there exists $\beta'_d$ such that
if $\beta > \beta'_d$ then the with probability going to $1$, the mixing time
of the dynamics on $G(n,d/n)$ is $\exp(\Omega(n))$.
\end{proposition}

\begin{proof}
The claim follows from expansion properties of $G(n,d/n)$. It is
well known that if $d > 1$ then with high probability $G(n,d/n)$
contains a {\em core} $C$ of size at least $\alpha_d n$ such that
that every $S \subset C$ of size at least $\alpha_d/4 n$ has at
least $\gamma_d n$ edges between $C$ and $S \setminus C$. Let $A$ be
the set of configurations $\sigma$ such that $\sigma$ restricted to
$C$ has at least $\alpha_d/4$ $+$'s and at least $\alpha_d/4$ $-$'s.
Then $P(A) \leq 2^n \exp(\beta |E| - 2 \beta \gamma_d n)/Z$. On the
other hand if $+$ denotes the all $+$ state then $P(+) = P(-) =
\exp(\beta |E|)/Z$. Thus by standard conductance arguments, the
mixing time is exponential in $n$ when $2 \exp(-2 \beta \gamma_d) <
1$.
\end{proof}

It is natural to conjecture that properties of the Ising model on
the branching process with $Poisson(d)$ offspring distribution
determines the mixing time of the dynamics on $G(n,d/n)$. In
particular, it is natural to conjecture that the critical point for
{\em uniqueness} of Gibbs measures plays a fundamental
role~\cite{Georgii:88,PemantleSteif:99} as results of similar flavor
were recently obtained for the hard-core model on random bi-partite
$d$ {\em regular} graphs~\cite{MoWeWo:07}. 

\begin{conjecture}
If $d \tanh(\beta) > 1$ then with high probability  over $G(n,d/n)$ the mixing time of the Gibbs sampler is
$\exp(\Omega(n))$. If $d > 1$ and $d \tanh(\beta) < 1$ then with high
probability over $G(n,d/n)$ the
mixing time of the Gibbs sampler is polynomial in $n$.
\end{conjecture}

After proposing the conjecture we have recently 
learned that Antoine Gerschenfeld and Andrea Montanari  
have found an elegant proof for estimating the partition function (that is the normalizing constant $Z(\beta)$) for the Ising model on random $d$-regular graphs~\cite{GerschenfeldMontanari:07}.  Their result together with a standard conductance argument shows exponentially slow mixing above the uniqueness threshold which in the context of random regular graphs is$(d+1)\tanh(\beta)=1$.

\subsection{Proof Technique}
Our proof follows the following main steps.
\begin{itemize}
\item
Analysis of the mixing time for Gibbs sampling on trees of varying degrees.
We find a bound on the mixing time on trees in terms of the maximal sum of
degrees along any simple path from the root. This implies that
for {\em all} $\beta$ if we consider
a tree where each node has number of descendants that has Poisson
distribution with parameter $d-1$ then with high probability the
mixing time of Gibbs sampling on the tree is polynomial in its size.
The motivation for this step is that we are looking at tree-like graphs
Note however, that the results established here hold for all
$\beta$, while rapid mixing for $G(n,d/n)$ does not hold for all
$\beta$. Our analysis here holds for all boundary conditions and all
external fields on the tree.
\item
We next use standard comparison arguments to extend the result above
to case where the graph is a tree with a few edges added. Note that
with high probability for all $v \in G(n,d/n)$ the induced subgraph
$B(v,\frac{1}{2} \log_d n)$ on all vertices of distance at most
$\frac{1}{2} \log_d n$ from $v$ is a tree with at most a few edges
added. (Note this still holds for all $\beta$).
\item
We next consider the effect of the boundary on the root
of the tree. We show that for tree of $a \log n$ levels,
the total variation distance of the conditional
distribution at the root given all $+$ boundary conditions and all
$-$ boundary conditions is $n^{-1-\Omega(1)}$ with probability
$1-n^{-1-\Omega(1)}$ provided $\beta < \beta_d$ is sufficiently
small (this is the only step where the fact that $\beta$ is small is
used).
\item
Using the construction of Weitz~\cite{Weitz:06} and a Lemma
from~\cite{KeMoPe:01,BeKeMoPe:05} we show that the spatial decay established
in the previous step also holds with probability $1-o(1)$ for all neighborhoods
$B(v,a \log  n)$ in the graph.
\item
The remaining steps use the fact that a strong enough decay of
correlation inside blocks each of which is rapidly mixing implies
that the dynamics on the full graph is rapidly mixing. This idea is
taken from~\cite{DSVW:04}.
\item
In order to show rapid mixing it suffices to exhibit a coupling of the dynamics
starting at all $+$ and all $-$ that couples with probability at least $1/2$ in
polynomial time. We show that the monotone coupling (where the configuration started at $-$ is always ``below'' the configuration started at $+$) satisfies this by showing that for each $v$ in polynomial time the two configurations at $v$
coupled except with probability $n^{-1}/(2e)$.
\item
In order to establish the later fact, it suffices to show that running
the dynamics on $B(v,a \log n)$ starting at all $+$ and all $+$
boundary conditions and the dynamics starting at all $-$ and all $-$ will couple at $v$ except with probability $n^{-1}/(2e)$ within polynomial time.
\item
The final fact then follows from the fact that the dynamics inside
$B(v,a \log n)$ have polynomial mixing time and that the stationary
distributions in $B(v,\frac{1}{2} \log_d n)$ given $+$ and $-$ boundary
conditions agree at $v$ with probability at least $1-n^{-1}/(4e)$.
\end{itemize}

We note that the decay of correlation on the self-avoiding tree defined
by Weitz that we prove here allows a different sampling scheme from the target
distribution. Indeed, this decay of correlation implies that given any
assignment to a subset of the vertices $S$ and any $v \not \in S$
we may calculate using the Weitz tree of radius $a \log n$ in polynomial time
the conditional probability that $\sigma(v) = +$ up to an
additive error of $n^{-1}/100$. It is easy to see that this allow sampling the
distribution in polynomial time. More specifically, consider the following 
algorithm from~\cite{Weitz:06}.

\begin{algorithm}\label{alg_weitz}
Fix a radius parameter $L$ and label the vertices $v_1,\ldots,v_n$.
Then the algorithm approximately samples from $P(\sigma)$ by
assigning the spins of $v_i$ sequentially.  Repeating from $1\leq i
\leq n$:
\begin{itemize}
\item In step $i$ construct $T_{SAW}^L(v_i)$, the tree of
self-avoiding walks truncated at distance $L$ from $v_i$.

\item Calculate
\[
p_i=P_{T_{SAW}^L}(\sigma_{v_i}=+|\sigma_{\{v_1,\ldots,v_{i-1}\}},\tau_{A-V_{i-1}}).
\]
(The boundary conditions at the tree can be chosen arbitrarily; in particular,
 one may calculate $p_i$ with no boundary conditions).

\item Fix $\sigma_{v_i}=X_{v_i}$ where $X_{v_i}$ is a random
variable with $p_i=P(X_{v_i}=+)=1-P(X_{v_i}=-)$.
\end{itemize}
\end{algorithm}

Then we prove that:

\begin{theorem}\label{thm_weitz_alg_random_graph}
Let $G$ be a random graph distributed as $G(n,d/n)$.  When
\[
\tanh(\beta)< \frac{1}{d},
\]
for any $\gamma>0$ there exist constants $r=r(d,\beta,\gamma)$ and
$C=C(d,\beta,\gamma)$ such that with high probability Algorithm
\ref{alg_weitz}, with parameter $r\log n$, has running time
$O(n^{C})$ and output distribution $Q$ with
$d_{TV}(P,Q)<n^{-\gamma}$. The result holds for the homogeneous
model~(\ref{eq:def_ising}) and for the inhomogeneous
model~(\ref{eq:def_ising_general}).
\end{theorem}

\begin{theorem}\label{thm_weitz_alg_general}
Let $G=(V,E)$ be any graph on $n$ vertices satisfying the following
properties. There exist $a >0, 0 < b < \infty$ such that for all $v
\in V$,
\begin{equation}\label{eq_SAW_growth}
|V_{T_{SAW}(v)}(v,a\log n)|\leq b^{a \log n}
\end{equation}
where $V_{T_{SAW}(v)}(v,r)=\{u\in T_{SAW}(v):d(u,v)\leq r\}$. When
\[
\tanh(\beta)< \frac{1}{b},
\]
for any $\gamma>0$ there exist constants $r=r(a,b,\beta,\gamma)$ and
$C=C(a,b,\beta,\gamma)$ such that Algorithm \ref{alg_weitz}, with
parameter $r\log n$, has running time $O(n^{C})$ and output
distribution $Q$ with $d_{TV}(P,Q)<n^{-\gamma}$. The result holds
for the homogeneous model~(\ref{eq:def_ising}) and for the
inhomogeneous model~(\ref{eq:def_ising_general}).
\end{theorem}

\subsection{Acknowledgment}
E.M. thanks Andrea Montanari and Alistair Sinclair
for interesting related discussions.  The authors thank Jinshan Zhang for pointing out an error in a previous draft.

\section{Proofs} \label{sec:proofs}

\subsection{Relaxation time on Sparse and Galton Watson Trees}
Recall that the local neighborhood of a vertex in $G(n,d/n)$ looks
like a branching process tree. In the first step of the proof we
bound the relaxation time on a tree generated by a Galton-Watson
branching process. More generally, we show that trees that are not
too dense have polynomial mixing time.

\begin{definition}
Let $T$ be a finite rooted tree.
We define $m(T) = \max_{\Gamma} \sum_{v \in \Gamma} d_v$ where the maximum is taken
over all simple paths $\Gamma$ emanating from the root and $d_v$ is the degree
of node $v$.
\end{definition}

\begin{theorem}\label{thm_relax_general_trees}
Let $\tau$ be the relaxation time of the continuous time Gibbs Sampler on $T$
where $0 \leq \beta_{u,v} \leq \beta$ for all $u$ and $v$ and given
arbitrary boundary conditions and external field. Then
\[
\tau \leq \exp(4 \beta m(T)).
\]
\end{theorem}

\begin{proof}

We proceed by induction on $m$ with a similar argument to the one
used in \cite{KeMoPe:01} for a regular tree. Note that if $m=0$ the claim
holds true since $\tau = 1$.
For the general case, let $v$ be the root of
$T$, and denote its children by $u_1, \ldots,u_k$ and denote the subtree
of the descendants of $u_i$ by $T^i$.  Now let $T'$ be the tree
obtained by removing the $k$ edges from $v$ to the $u_i$, let $P'$
be the Ising model on $T'$ and let $\tau'$ be the relaxation time
on $T' $. By equation \eqref{eq_relax_defn} we have that
\begin{equation} \label{eq:recursion_ising}
\tau /\tau' \leq \frac{\max_\sigma
P(\sigma)/P'(\sigma)}{\min_{\sigma,\tau}Q(\sigma,\tau)/Q'(\sigma,\tau)}
\leq \exp(4 \beta k).
\end{equation}
Now we divide $T'$ into $k+1$
blocks $\{\{v\},\{T^1\},\ldots,\{T^k\} \}$. Since
these blocks are not connected to each other the block dynamics is simply the product chain.  Each block updates at rate 1 and therefore the relaxation time of the block dynamics is simply $1$.
By applying Proposition\ref{prop:martinelli} we get that the relaxation time on $T'$ is
simply the maximum of the relaxation times on the blocks,
\[
\tau' \leq \max \{1,\tau^i\}.
\]
where $\tau^i$ is the relaxation time on $T^i$.
Note that by the definition of $m$, it follows that the value of $m$ for each
of the subtrees $T^i$ satisfies $m(T^i) \leq m - k$, and therefore for all
$i$ it holds that $\tau^i \leq \exp(4 \beta (m-k))$.
This then implies by~(\ref{eq:recursion_ising}) that
$\tau \leq \exp(4 \beta m)$ as needed.
\end{proof}

\subsection{Some properties of Galton Watson Trees}
Here we prove a couple of useful properties for Galton Watson trees that will
be used below. We let $T$ be the tree generated by a Galton-Watson branching process
with offspring distribution $N$ such that for all $t$, $E\exp(t N) <
\infty$ and such that $E(N) = d$. Of particular interest to us would be 
the Poisson distribution with mean $d$ which has
\[
E\exp(t N) = \exp(d (e^t-1)).
\]
We let $T_r$ denote the first $r$ levels of $T$. We let $M(r)$ denote
the value of $m$ for $T(r)$ and $\tau(r)$ the supremum of the relaxation times of the continuous time Gibbs Sampler on $T(r)$ over any boundary conditions and external fields assuming
that $\beta = \sup \beta_{u,v}$. We denote by $Z_r$
the number of descendants at level $r$.

\begin{theorem}\label{thm_relax_time_GW_tree}
Under the assumptions above we have:
\begin{itemize}
\item
There exists a positive function $c(t)$ such that for all $t$ and all $r$:
\[
E[\exp(t M(r))] \leq \exp(c(t) r).
\]
\item
Then $E \tau(r)\leq {C(\beta)}^r$
for some $C(\beta)<\infty$ depending on $\beta = \sup \beta_{u,v}$
only.
\item If $N$ is the Poisson distribution with mean $d$ then for all $t>0$,
\[
\sup_r E[\exp(t Z_r d^{-r})] < \infty.
\]
\end{itemize}
\end{theorem}

\begin{proof}
Let $K$ denote the degree of the root of $T_r$ and for $1 \leq i
\leq K$ let $M_i(r-1)$ denote the value of $m$ for the sub-tree of
$T_r$ rooted at the $i$'th child. Then:
\begin{eqnarray*}
E[\exp( t M(r))] &=& E[\max(1,\max_{1 \leq i \leq K} \exp(t
(M_i(r-1) + K)))] \\ &\leq& E[(1+\exp(t K)) \sum_{i=1}^K \exp(t
M_i(r-1))]  \\ &=& E[(1+K \exp(t K))] E[\exp(t M(r-1))].
\end{eqnarray*}
and so the result follows by induction provided that $c(t)$ is large
enough so that 
\[
\exp(c(t)) \geq E(1+K \exp(t K)).
\]

\bigskip

For the second statement of the theorem, note that
by the previous theorem we have that 
\[
E \tau(r) \leq E[\exp( 4 \beta M(r))], 
\]
where $M(r)$ is the random value of $m$ for the tree $T_r$
so if $C(\beta)=\exp(c(4\beta))$ then $E \tau(r)\leq {C(\beta)}^r$.

\bigskip

For the last part of the theorem, let $N_i$ be independent copies of $N$ 
and note that
\begin{eqnarray}\label{eq_exp_moment_recursion}
E\exp(t Z_{r+1}) &=& E\exp(\sum_{i=0}^{Z_r} t d^{-(r+1)} N_i) 
= E[E[\exp(\sum_{i=0}^{Z_r} t d^{-(r+1)} N_i | Z_n]] \\ \nonumber
&=& 
E[(E[\exp(t d^{-r+1} N)])^{Z_r}] =   
E \exp(\log(E\exp(t d^{-(r+1)} N))Z_r)
\end{eqnarray}
which recursively relates the exponential moments of $Z_{r+1}$ to
the exponential moments of $Z_r$.  In particular since all the
exponential moments of $Z_1$ exist, $E \exp(t Z_{r})<\infty$ for all
$t$ and $r$. When $0 < s \leq 1$
\begin{equation}\label{eq_exp_moment_bound}
E\exp(s N) = \sum_{i=0}^\infty \frac{s^i E N^i}{i!} \leq 1 + sd +
s^2 \sum_{i=2}^\infty \frac{E N^i}{i!} \leq \exp(sd(1+\alpha s))
\end{equation}
provided $\alpha$ is sufficiently large.  Now fix a $t$ and let
$t_n=t\exp(2\alpha t \sum_{i=r+1}^\infty d^{-i})$. For some
sufficiently large $j$ we have that $\exp(2\alpha t
\sum_{i=r+1}^\infty d^{-i})<2$ and $t_r d^{-(r+1)}<1$ for all $r\geq
j$. Then for $r\geq j$ by equations~\eqref{eq_exp_moment_recursion}
and \eqref{eq_exp_moment_bound},
\begin{align*}
E\exp(t_{r+1} Z_{r+1} d^{-(r+1)}) &= E \exp(\log(E\exp(t_{r+1}
d^{-(r+1)} N_i))Z_r)\\
&\leq E\exp(t_{r+1}(1+\alpha t_{r+1} d^{-(r+1)}) Z_{r} d^{-r})\\
&\leq E\exp(t_{r+1}(1+2\alpha t d^{-(r+1)}) Z_{r} d^{-r})\\
&\leq E\exp(t_{r} Z_{r} d^{-r})
\end{align*}
and so
\[
\sup_{r\geq j} E\exp(t Z_{r} d^{-r}) \leq \sup_{r\geq j} E\exp(t_{r}
Z_{r} d^{-r}) = E\exp(t_{j} Z_{j} d^{-j}) < \infty
\]
which completes the result.
\end{proof}

When the branching process is super-critical,
the number of vertices is $O((EW)^r)$ and the result above gives that the
mixing time is polynomial in the number of vertices on Galton Watson
branching process with high probability.  We remark that all our
bounds here are increasing in the degrees of the vertices so if a
random tree T is stochastically dominated by a Galton-Watson
branching process then the same bound applies.

\subsection{Relaxation in Tree-Like Graphs}

For the applications considered for random and sparse graphs,
it is not always the case that the neighborhood of a vertex is a tree, instead
it is sometimes a tree with a small number of edges added.
Using standard comparison arguments we show that the mixing time of a graph that is a tree with a few edges added is still polynomial. We also show that
with high probability for the $G(n,d/n)$ the neighborhoods of all vertices
are tree-like.

\begin{proposition}\label{spanning_tree}
Let $G$ be a graph on $r$ vertices with $r+s-1$ edges that has a
spanning tree $T$ with $m(T) = m$. Then the mixing time $\tau$ of
the Glauber dynamics on $G$ with any boundary conditions and
external fields satisfies:
\[
\tau \leq \exp(4 \beta(m+s)).
\]
\end{proposition}

\begin{proof}
By equation \eqref{eq:recursion_ising} removing the $s$ edges in $G$ which are not in $T$ decreases the relaxation time by at most a multiplicative factor of $\exp(4 \beta s)$.  By Theorem \ref{thm_relax_general_trees} the relaxation time of $T$ is at most $\exp(4 \beta m)$ so the relaxation time of $G$ is bounded by $\exp(4 \beta(m+s))$.
\end{proof}

\begin{lemma}\label{lem_random_graph_to_general}
Let $G$ be a random graph distributed as $G(n,d/n)$.  The following hold with high probability over $G$:
\begin{itemize}
\item For $0<a<\frac1{2\log d}$ there exists some $c(a,d)$ such that for all $v\in G$, $m(v,a \log n)\leq c \log n$. 

\item There
exists $k=k(a,d)>0$ such that for all $v\in
G$, $t(v,a \log n)\leq k$.
\item For $0<a<\frac1{2\log d}$ and every $v\in G$,
\[
|B(v,a\log n)|\leq
3(1-d^{-1})n^{a\log d} \log n.
\]
\end{itemize}
\end{lemma}

\begin{proof}
We construct a spanning tree $T(v,l)$ of $B(v,l)$ in a standard
manner. Take some arbitrary ordering of the vertices of $G$.  Start
with the vertex $v$ and attach it to all its neighbors in $G$.  Now
take the minimal vertex in $S(v,1)$, according to the ordering, and
attach it to all its neighbors in G which are not already in the
graph.  Repeat this for each of the vertices in $S(v,1)$ in
increasing order.  Repeat this for $S(v,2)$ and continue until
$S(v,l-1)$ which completes $T(v,l)$. By construction this is a
spanning tree for $B(v,l)$. The construction can be viewed as a
breadth first search of $B(v,l)$ starting from $v$ and exploring
according to our ordering.

By a standard argument $T(v,a \log n)$ is stochastically dominated
by a Galton-Watson branching process with offspring distribution
Poisson$(d)$.  Then by repeating the argument of Theorem
\ref{thm_relax_time_GW_tree} for some $\delta$,
\[
E\exp (m(T(v,a \log
n),v,a \log n))\leq \delta^{a\log n}
\]
 and so,
\[
P(m(T(v,a \log n),v,a \log n))>(a\delta+2)\log n)=O(n^{-2}).
\]
which implies that with high probability $m(T(v,a \log n),v,a \log
n))<(a\delta+2)\log n$ for all $v$.

If $Z_l$ are the number of offspring in generation $l$ of a
Galton-Watson branching process with offspring distribution
Poisson$(d)$ then by Theorem~\ref{thm_relax_time_GW_tree} we have that
$\sup_l E \exp(Z_l/d^l)< \infty$ and since
\[
P(|S(v,l)| > 3d^l \log n)\leq  P(\exp(Z_l/d^l) > n^3)  \leq n^{-3} E
\exp(Z_l/d^l),
\]
it follows by a union bound over all $v\in G$ and $1\leq l \leq a\log n$
we have with high probability for all $v$, 
\begin{equation}\label{e:ERgraphGrowthBound}
|B(v,a\log n)|\leq
3(1-d^{-1})n^{a\log d} \log n.
\end{equation}

In the construction of $T(v,a\log n)$ there may be some edges in
$B(v,a\log n)$ which are not explored and so are not in $T(v,a\log
n)$.  Each edge between $u,w\in V(v,a\log n)$ which is not explored
in the construction of $T(v,a\log n)$ then is present in $B(v,a\log
n)$ independently of $T(v,a\log n)$ with probability $d/n$.  There
are at most $(3(1-d^{-1})n^{a\log d} \log n)^2$ unexplored edges.
Now when $k>1/(1-2a\log d)$,
\[
P(\hbox{Binomial}((3(1-d^{-1})n^{a\log d} \log n)^2,d/n)>k)=
O(n^{k(2a\log d-1)} (\log n)^{2k})=n^{-1-\Omega(1)}
\]
so by a union bound with high probability we have $t(v,a \log n)\leq
k$.  Now a self-avoiding path in $B(v,a\log n)$ can traverse each of
these $k$ edges at most once so this path can be split into at most
$k+1$ self-avoiding paths in $T(v,a\log n)$ and hence with high
probability $m(v,l)\leq c \log n$ where $c=(k+1)(a\delta+2)$.
\end{proof}

\begin{lemma}\label{lem_SAW_growth_random_graph}
When $0<a<\frac1{2\log d}$ with high probability for all $v\in G$,
\[
|V_{T_{SAW}(v)}(v,a\log n)|\leq O(n^{a\log d} \log n)
\]
where $V_{T_{SAW}(v)}(v,r)=\{u\in T_{SAW}(v):d(u,v)\leq r\}$.
\end{lemma}

\begin{proof}
We now count the number of self-avoiding walks of length at most
$a\log n$ in $B(v,a\log n)$.  By
Lemma~\ref{lem_random_graph_to_general} we have that with high
probability for all $v$, $|B(v,a\log n)|\leq 3(1-d^{-1})n^{a\log d}
\log n$ and $t(v,a \log n)\leq k$.  Let $e_1,\ldots,e_{t(v,a \log
n)}$ denote the edges in $B(v,a\log n)$ which are not in $T(v,a\log
n)$. Now every vertex in $u'\in T_{SAW}(v)$ corresponds to a unique
self avoiding walk in $B(v,a\log n)$ from $v$ to $u$.  A
self-avoiding walk in $B(v,a\log n)$ passes through each edge at
most most once so in particular it passes through each of the $e_i$
at most once.  So a path which begins at $v$ traverses through some
sequence $e_{i_1},\ldots,e_{i_l}$ in particular directions and then
ends at $u$ is otherwise uniquely defined since the intermediate
steps are paths in $T(v,a\log n)$ which are unique.  There are at
most $k(k!)$ sequences $e_{i_1},\ldots,e_{i_l}$, there are $2^k$
choices of directions to travel through them, and at most
$3(1-d^{-1})n^{a\log d} \log n$ possible terminal vertices in
$B(v,a\log n)$ so $|V_{T_{SAW}(v)}(v,a\log n)|\leq 3(1-d^{-1})2^k
k(k!)n^{a\log d} \log n$.

\end{proof}

\subsection{Spatial decay of correlation for tree-like neighborhoods}

\begin{proposition}\label{prop_tree_size}
Let $T$ be a tree such that $m(v,a )\leq m$.  Then $S|(v,a)|\leq
(\frac{m-a+1}{a})^{a}$.
\end{proposition}

\begin{proof}
First we establish inductively that $|S(v,a)|$ is
maximized by a spherically symmetric tree, that is one where the
degrees of the vertices depend only on their distance to $v$ (it may be that it is also maximized by non-spherically symmetric trees).  It is clearly true when $a=0$ so suppose that it is true for all $m$ up to height $a-1$.  Let $T^*$ be a tree of height $a$ rooted at $v$ that maximizes $|S(T^*,v,a)|$ under the constraint $m(T^*,v,a )\leq m$ and let $k$ be the degree of $v$.  Then each of the subtrees $T_i$ attached to $v$ have depth $a-1$ and are constrained to have $m(T_i,v_i,a -1)\leq m-k-1$.
Let $T^{-}$ be a sphereically  symmetric tree of height $a-1$ which has $m(T^{-},v,a -1)\leq m-k-1$ and maximizes $S|(T^{-},v,a-1)|$.  A vertex $v$ connected to the roots of $k$ copies of $T^{-}$ is a spherically symmetric tree of height $a$ with $m(v,a)=m$ and by our inductive hypothesis must have boundary size $|S(v,a)|$ at least as large as $T^*$ which completes the induction step.

So suppose that $T$ is sphereically symmetric and let
$d_i$ be the degree of a vertex distance $i$ from $v$.  Then by the
arithmetic-geometric inequality
\[
|S(v,a)|=d_0\prod_{i=1}^{a-1} (d_i-1) \leq ((\sum_{i=0}^{a-1} d_i
-(a-1))/a)^{a} \leq \left(\frac{m-a+1}{a}\right)^{a}.
\]
\end{proof}

We now consider the effect that conditioning on the leaves of a tree can have on the marginal distribution of the spin at the root.  It will be convenient to compare this probability to the Ising model with the same interaction strengths $\beta_{uv}$ but no external field ($h\equiv 0$) which we will denote $\widetilde{P}$.

\begin{lemma}\label{lem_tree_boundary_effect}
If $T$ is a tree, $P$ is the Ising model with arbitrary external field
(including $h_u = \pm \infty$ meaning that $\sigma_u$ is set to  $\pm$) and
$\beta_{u,v} \leq \beta$ then for all $v$ ,
\[
P(\sigma_v=+|\sigma_{S(v,l)}\equiv +)-
P(\sigma_v=+|\sigma_{S(v,l)}\equiv -) \leq |S(v,l)|  (\tanh
\beta)^l.
\]
\end{lemma}

\begin{proof}
Lemma 4.1
of~\cite{BeKeMoPe:05} states that for any vertices $v,u\in T$,
\begin{equation}\label{e:noFieldEffect}
P(\sigma_v=+|\sigma_{u}=+)-
P(\sigma_v=+|\sigma_{u}=-) \leq \widetilde{P}(\sigma_v=+|\sigma_{u}=+)-
\widetilde{P}(\sigma_v=+|\sigma_{u}=-).
\end{equation}
If $u_0,u_1,\ldots,u_l$ are a path of vertices in $T$ then a simple calculation yields that
\begin{equation}\label{e:pathDecay}
\widetilde{P}(\sigma_{u_k}=+|\sigma_{u_0}=+)-
\widetilde{P}(\sigma_{u_k}=+|\sigma_{u_0}=-) = \prod_{i=1}^k \tanh \beta_{u_{i-1} u_i} \leq (\tanh \beta)^k.
\end{equation}
Now suppose that $u\in S(v,l)$ and that $\eta_{S(v,l)}^+$ and $\eta_{S(v,l)}^-$ are configurations on $S(v,l)$ which differ only at $u$ where $\eta_u^\pm = \pm$.  Conditioning is equivalent to setting an infinite external field so equations \eqref{e:noFieldEffect} and \eqref{e:pathDecay}  imply that
\begin{equation}\label{e:singleSpinFlip}
P(\sigma_v=+|\sigma_{S(v,l)}=\eta^+)- P(\sigma_v=+| \sigma_{S(v,l)}=\eta^-) \leq (\tanh \beta)^l.
\end{equation}
Take a sequence of configurations $\eta^0,  \eta^1,\ldots, \eta^{|S(v,l)|}$ on $S(v,l)$ with $\eta^0\equiv -$ and $ \eta^{|S(v,l)|} \equiv+$  where consecutive configurations differ at a single vertex.  By equation \eqref{e:singleSpinFlip} we have that
\[
P(\sigma_v=+|\sigma_{S(v,l)}= \eta^{i+1})-
P(\sigma_v=+|\sigma_{S(v,l)}= \eta^i) \leq (\tanh
\beta)^l
\]
and so
\[
P(\sigma_v=+|\sigma_{S(v,l)}\equiv +)-
P(\sigma_v=+|\sigma_{S(v,l)}\equiv -) \leq |S(v,l)|  (\tanh
\beta)^l
\]
which completes the proof.
\end{proof}

Now $B(v,a\log n)$ is not in general a tree so we use the
self-avoiding tree construction of Weitz \cite{Weitz:06} to reduce
the problem to one on a tree.  The tree of self-avoiding walks,
which we denote $T_{saw}(v,a\log n)$, is the tree of paths in
$B(v,a\log n)$ starting from $v$ and not intersecting
themselves, except at the terminal vertex of the path. Through this
construction each vertex in $T_{saw}(v,a\log n)$ can be identified
with a vertex in $G$ which gives a natural way to relate a subset
$\Lambda\subset V$ and a configuration $\sigma_\Lambda$ to the
corresponding subset $\Lambda'\subset T_{saw}(v,a\log n)$ and configuration
$\sigma_{\Lambda'}$ in $T_{saw}$. Furthermore if $A,B\subset V$ then
$d(A,B)=d(A',B')$. Then Theorem 3.1 of \cite{Weitz:06} gives the
following result.  Each vertex (edge) of $T_{saw}$ corresponds to a
vertex (edge) so $P_{T_{saw}}$ is defined by taking the
corresponding external field and interactions.

\begin{lemma}\label{lem_Weitz}
For a graph $G$ and $v\in G$ there exists $A\subset T_{saw}$ and
some configuration $\tau_A$ on $A$ such that,
\[
P_G(\sigma_v=+|\sigma_\Lambda)=P_{T_{saw}}(\sigma_v=+|\sigma_{\Lambda'},\tau_{A-\Lambda'}).
\]
The set $A$ corresponds to the terminal vertices of path which
returns to a vertex already visited by the path.
\end{lemma}

\begin{corollary}\label{coroll_ball_boundary_effect}
Suppose that $a,b,c,\beta$ satisfy the hypothesis of Theorem
\ref{thm_mix_general_graph}. Then,
\[
\max_{v\in G} P(\sigma_v=+|\sigma_{S(v,a \log n)}\equiv +)-
P(\sigma_v=+|\sigma_{S(v,a \log n)}\equiv -) = o(n^{-1}).
\]
\end{corollary}

\begin{proof}
By applying Lemma \ref{lem_Weitz} we have that if $\Lambda=S(v,a \log n)$
then
\begin{align*}
&P_G(\sigma_v=+|\sigma_\Lambda\equiv
+) - P_G(\sigma_v=+|\sigma_\Lambda\equiv
-)\\
&\qquad=P_{T_{saw}}(\sigma_v=+|\sigma_{\Lambda'}\equiv
+,\tau_{A-\Lambda'}) - P_{T_{saw}}(\sigma_v=+|\sigma_{\Lambda'}\equiv
-,\tau_{A-\Lambda'}).
\end{align*}
Conditioning on $\tau_A$ is equivalent to setting the external field
at $u\in A$ to $\hbox{sign}(\tau_v)\infty$ hence it follows by Lemma \ref{lem_tree_boundary_effect} that,
\[
P_{T_{saw}}(\sigma_v=+|\sigma_{\Lambda'}\equiv
+,\tau_{A-\Lambda'}) - P_{T_{saw}}(\sigma_v=+|\sigma_{\Lambda'}\equiv
-,\tau_{A-\Lambda'})  \leq |S_{saw}(v,a\log n)| (\tanh \beta)^{a\log n}
\]
where $S_{saw}(v,a\log n)=\{u\in T_{saw}(v,a\log n):d(u,v)=a\log
n\}$. Now suppose $v=u_1,u_2,\ldots,u_{k}$ is a non-repeating walk
in $T_{saw}$ and let $u'_1,u'_2,\ldots,u'_{k}$ be the corresponding
walk in $G$. Then from the construction of $T_{saw}$ either
$u'_1,u'_2,\ldots,u'_{k}$ is a non-repeating walk in $G$ or for some
$j<k$, $u'_j=u'_k$ in which case $u_k$ is a leaf of $T_{saw}$ and so
has degree 1.  It also follows from the construction of $T_{saw}$
that the degree of $u_i$ is less than or equal to the degree of
$u'_i$ and so we have that $m(v,a\log n)\leq m(T_{saw},v,a\log
n)+1$. The by Proposition \ref{prop_tree_size} $$|S_{saw}(v,a\log
n)|\leq \left(\frac{c\log n-a\log n +2}{a\log n}\right)^{a\log n}=O\left(n^{a\log
((c-a)/a)}\right)=o(n^{-1}(\tanh \beta)^{-a\log n}),$$  which completes the result.
\end{proof}

\subsection{Proof of the Main Result}

\begin{proof}(Theorem \ref{thm_mix_general_graph})
Let $X_t^{+},X_t^{-}$, denote the Gibbs sampler on $G$ started from respectively all $+$ and $-$,
coupled  using the monotone coupling described in Section \ref{s:monotoneCoupling}. Fix some vertex $v\in
G$. Define four new chains $Q_t^+$, $Q_t^-$, $Z_t^+$ and $Z_t^-$.
These chains run the Glauber dynamics and are coupled with $X_t^+$
and $X_t^-$ inside $B(v,a \log n)$ by using the same choice of vertices $v_1,v_2,\ldots$ and the same choice of update random variables $U_1,U_2,\ldots$ except that they are fixed (i.e. do not
update) outside $B(v,a \log n)$. They are given the following
initial and boundary conditions.
\begin{itemize}

\item $Q_t^+$ starts from all $+$ configuration (and therefore has all
$+$ boundary conditions during the dynamics).

\item $Q_t^-$ starts from all $-$ configuration (and therefore has all
$-$ boundary conditions during the dynamics).

\item $Z_t^+$ starts from all $+$ configuration outside $B(v,a \log
n)$ and $Z_0^+$ is distributed according to the stationary
distribution inside $B(v,a \log n)$ given the all $+$ boundary
condition (therefore $Z_t^+$ will have this distribution for all
$t$).

\item $Z_t^-$ starts from all $-$ configuration outside $B(v,a \log
n)$ and is distributed according to the stationary distribution
inside $B(v,a \log n)$ given the all $-$ boundary condition
(therefore $Z_t^-$ will have this distribution for all $t$).
\end{itemize}
As the Gibbs distribution on $B(v,a \log n)$ with a $+$ boundary condition stochstically dominates the distribution with a $-$ boundary condition, we can initialize $Z_t^+$ and $Z_t^-$ so that $Z_0^+ \succcurlyeq
Z_0^-$.  By  monotonicity of the updates we have $Q_t^+ \succcurlyeq Z_t^+ \succcurlyeq Z_t^- \succcurlyeq Q_t^- $ for
all $t$.  We also have that $Q_t^+ \succcurlyeq X_t^+ \succcurlyeq
X_t^- \succcurlyeq Q_t^-$ on  $B(v,a \log n)$.  As $Z_t^+$ (respectively $Z_t^-$) starts in the stationary distribution of the Gibbs sampler given the all $+$ (respectively all $-$) boundary condition, it remains in the stationary distribution for all time $t$.

Since $Z_t^+(v) \geq Z_t^-(v)$ we have that
\begin{align*}
P(Z_t^+(v) \neq Z_t^-(v)) = P(Z_t^+(v) = + ) -P( Z_t^-(v)=+) \leq o(n^{-1}),
\end{align*}
for all $t$ where the inequality follows from Corollary~\ref{coroll_ball_boundary_effect}.  By Proposition \ref{spanning_tree}
the continuous time Gibbs sampler on $B(v,a\log n)$ has relaxation time bounded above by $\exp(4\beta(b+c)\log n)$ which implies that  the discrete time relaxation time satisfies $\tau \leq n^{1+4\beta(b+c)}$.
As each vertex has degree at most $c\log n$,
\[
\log(\min_\sigma P(\sigma))^{-1}\leq (\beta|E|) + \sum_u |h_u|
\leq (100 c n^2 \beta^2 \log n)
\]
which implies that $\tau_{mix}\leq O(n^{4+4(b+c)\beta})$ since the mixing satisifies $\tau_{mix}\leq
\tau(1+\frac12 \log(\min_\sigma P(\sigma))^{-1})$. 
For $C = 6+4(b+c)\beta$ we have that
with high probability after $t=2 n^C$ steps that the Gibbs sampler
has chosen every vertex at least $n^{5+4(b+c)\beta}\geq n\tau_{mix}$ times.  It follows that the number of updates to $B(v,a\log n)$ is at least $n$ times its mixing time and so
\[
P(Q_t^+(v)\neq Z_t^+(v)) \leq d_{TV}(Q_t^+(v) , Z_t^+(v)) \leq e^{-n} = o(n^{-1}).
\]
where $d_{TV}$ denotes the total variation distance which is always bounded above by $\exp(-t/\tau_{mix})$.  We similarly have that
\[
P(Q_t^-(v)\neq Z_t^-(v))\leq o(n^{-1}).
\]
It follows that $P(Q_t^+(v)\neq Q_t^-(v))\leq o(n^{-1})$ and hence
$P(X_t^+(v)\neq X_t^-(v))\leq o(n^{-1})$ for all $v$. By a union bound
$P(X_t^+\neq X_t^-)\leq o(1)$ so the mixing time is bounded by
$O(n^C)$ as required.
\end{proof}

\begin{proof}(Theorem \ref{thm_mix_time_random_graph})
By Lemma~\ref{lem_random_graph_to_general} with high probability a
random graph satisfies the hypothesis of
Theorem~\ref{thm_mix_general_graph} for small enough $\beta$. To
prove the result when $\tanh(\beta)< \frac{1}{e^2 d}$ the only
modification to the proof of Theorem \ref{thm_mix_general_graph}
needed is to show that with high probability when $-1/(\log (d\tanh
(\beta))) < a < (2 \log d)^{-1}$ we still have $P(Z_t^+(v) \neq
Z_t^-(v))\leq o(n^{-1})$.  We know from
Lemma~\ref{lem_SAW_growth_random_graph} that with high probability
$|V_{T_{SAW}(v)}(v,a\log n)|\leq O(n^{a\log d} \log
n)=o(n^{-1}(\tanh \beta)^{-a\log n})$. Now using this bound and
repeating the proof of Corollary~\ref{coroll_ball_boundary_effect}
we get that $P(Z_t^+(v) \neq Z_t^-(v))= o(n^{-1})$ as required.

The mixing time is bounded by $n^{6+4(b+c)\beta}$ which is bounded by $n^{6+4(b+c)\tanh^{-1}(\frac1{e^2 d})}$ and does not need to depend on $\beta$.
\end{proof}

\subsection{Sampling from the distribution through the tree of self avoiding walks}

The proofs Theorems~\ref{thm_weitz_alg_random_graph} 
and~\ref{thm_weitz_alg_general} make use the following lemmas.
\begin{lemma}\label{lemma_TV_distance}
Let $(X_1,\ldots,X_n)$ and $(Y_1,\ldots,Y_n)$
be two vector valued distributions taking values in some product
space. Suppose that for all $1 \leq i \leq n$ and all $(x_1,\ldots,x_i)$
we have
\[
d_{TV}((X_i | X_1 = x_1,\ldots, X_{i-1} = x_{i-1}),
      (Y_i | Y_1 = x_1,\ldots, Y_{i-1} = x_{i-1})) \leq \eps_i,
\]
Then
\[
d_{TV}((X_1,\ldots,X_n),(Y_1,\ldots,Y_n)) \leq \sum_{i=1}^n \eps_i,
\]
\end{lemma}
\begin{proof}
The proof follows by constructing a coupling of the two
distributions whose total variation distance is bounded by
$\sum_{i=1}^n \eps_i$. The coupling is performed by first coupling
$X_1$ and $Y_1$ except with probability $\eps_1$. Then at step $i$,
given the coupling of $(X_1,\ldots,X_{i-1})$ and
$(Y_1,\ldots,Y_{i-1})$ and conditioned on
\[
(X_1,\ldots,X_{i-1}) = (Y_1,\ldots,Y_{i-1}),
\]
we couple the two configurations in such a way that they do not agree at
most with probability $\eps_i$. The proof follows.
\end{proof}

\begin{lemma}\label{Big_SAW_Tree}
Suppose the graph $G$ satisfies that for $v \in V$ it holds that
\[
|V_{T_{SAW}(v)}(v,a)|\leq b,
\]
Then for all integer $j$ it holds that
\[
|V_{T_{SAW}(v)}(v,ja)|\leq b^j.
\]
\end{lemma}

\begin{proof}
We prove the result by induction on $j$. Suppose that $u\in
S_{T_{SAW}(v)}(v,(j-1)a)$ and let $T_u$ denote the subtree of $u$
and its descendants in $V_{T_{SAW}(v)}(v,j a) \backslash
V_{T_{SAW}(v)}(v,(j-1) a)$.  Each path from $u$ in $T_u$ corresponds
to a self avoiding walk in $G$ started from $u$ so it follows that
the number of vertices in $T_u \backslash \{u\}$ is at most $b-1$.
So $|V_{T_{SAW}(v)}(v,j a) \backslash V_{T_{SAW}(v)}(v,(j-1) a)|\leq
b^{j-1}(b-1)$ which completes the induction.
\end{proof}



\begin{proof}(Theorem \ref{thm_weitz_alg_general})
Set $r=j a$ where $j$ is the smallest integer greater than
$\frac{-(1+\gamma)}{\log(b\tanh \beta)}$.  By Lemma
\ref{Big_SAW_Tree} for all $i$, $|V_{T_{SAW}(v_i)}(v_i,r \log
n)|\leq b^{r\log n}$ so $T_{SAW}^{r\log n}(v_i)$ the tree of self
avoiding walks of radius $r\log n$ can be constructed in $O(b^{r\log
n})=O(n^{r\log b})$ steps. Using the standard recursions on a tree,
$p_i$ can be evaluated in $O(n^{r\log b})$ steps so the running time
of the algorithm is $O(n^C)$ where $C=1+r\log b$.

At step $i$ we calculate $p_i=P_{T_{SAW}^{r\log
n}}(\sigma_{v_i}=+|\sigma_{V_{i-1}},\tau_{A-V_{i-1}})$ to
approximate $P(\sigma_{v_i}=+|\sigma_{V_{i-1}})$.  Applying Lemma
\ref{lem_Weitz} we have that
\[
P(\sigma_{v_i}=+|\sigma_{V_{i-1}}) = P_{T_{SAW}(v_i)}
(\sigma_{v_i}=+| \sigma_{V_{i-1}},\tau_{A-V_{i-1}}).
\]
where $V_j=\{v_1,\ldots,v_{j}\}$ and so if
$\Lambda=S_{T_{SAW}(v_i)}(v_i,r\log n)$ then,
\begin{align*}
P_{T_{SAW}(v_i)} (\sigma_{v_i}=+| \sigma_{\Lambda}\equiv -,
\sigma_{V_{i-1}}, \tau_{A-V_{i-1}}) &\leq
P(\sigma_{v_i}=+|\sigma_{V_{i-1}}) \\
&\leq P_{T_{SAW}(v_i)} (\sigma_{v_i}=+| \sigma_{\Lambda}\equiv +,
\sigma_{V_{i-1}}, \tau_{A-V_{i-1}})
\end{align*}
and similarly
\begin{align*}
P_{T_{SAW}(v_i)} (\sigma_{v_i}=+| \sigma_{\Lambda}\equiv -,
\sigma_{V_{i-1}}, \tau_{A-V_{i-1}}) &\leq P_{T_{SAW}^{r\log
n}}(\sigma_{v_i}=+|\sigma_{V_{i-1}},\tau_{A-V_{i-1}})\\
&\leq P_{T_{SAW}(v_i)} (\sigma_{v_i}=+| \sigma_{\Lambda}\equiv +,
\sigma_{V_{i-1}}, \tau_{A-V_{i-1}})
\end{align*}
so
\begin{align*}
&|P_{T_{SAW}^{r\log
n}}(\sigma_{v_i}=+|\sigma_{V_{i-1}},\tau_{A-V_{i-1}})
- P(\sigma_{v_i}=+|\sigma_{V_{i-1}})|\\
&\leq  P_{T_{SAW}(v_i)} (\sigma_{v_i}=+| \sigma_{\Lambda}\equiv +,
\sigma_{V_{i-1}}, \tau_{A-V_{i-1}}) - P_{T_{SAW}(v_i)}
(\sigma_{v_i}=+| \sigma_{\Lambda}\equiv -, \sigma_{V_{i-1}},
\tau_{A-V_{i-1}}).
\end{align*}
Conditioning on $\sigma_{V_{i-1}}$ and $\tau_A$ is equivalent to
setting the external field to be $\pm\infty$.  Then by Lemma \ref{lem_tree_boundary_effect}
\begin{align*}
&P_{T_{SAW}(v_i)} (\sigma_{v_i}=+| \sigma_{\Lambda}\equiv +,
\sigma_{V_{i-1}}, \tau_{A-V_{i-1}})  -  P_{T_{SAW}(v_i)}
(\sigma_{v_i}=+| \sigma_{\Lambda}\equiv -, \sigma_{V_{i-1}},
\tau_{A-V_{i-1}})  \\
&\qquad \leq |S_{T_{SAW}(v_i)}(v_i,r\log n)| (\tanh \beta)^{r \log n} =  O(n^{-1-\gamma}).
\end{align*}
If $Q$ is
the output of the algorithm then by Lemma \ref{lemma_TV_distance}
\[
d_{TV}(P,Q)\leq \sum_{i=1}^n
\sup_{\sigma_{V_{i-1}}}|P_{T_{SAW}^{r\log
n}}(\sigma_{v_i}=+|\sigma_{V_{i-1}},\tau_{A-V_{i-1}}) -
P(\sigma_{v_i}=+|\sigma_{V_{i-1}})|= O(n^{-\gamma})
\]
which completes the result.
\end{proof}

\begin{proof}(Theorem~\ref{thm_weitz_alg_random_graph})
By Lemma~\ref{lem_SAW_growth_random_graph}
equation~\eqref{eq_SAW_growth} holds with high probability for any
$0<a<\frac1{2\log d}$ and $b>d$ so the result follows by
Theorem~\ref{thm_weitz_alg_general}.

\end{proof}

\clearpage

\bibliographystyle{plain}
\bibliography{all,my}

\end{document}